\documentclass[12pt]{article}
\usepackage{epsfig}
\begin{document}

\def\Empty{}
%
%
\def\bottomfraction{.7}
\def\textfraction{0}
\def\floatpagefraction{.7}
%
\catcode`\@=11

\def\section{\@startsection {section}{1}{\z@}{-3.5ex plus-1ex minus
    -.2ex}{2.3ex plus.2ex}{\large\bf}}
\def\subsection{\@startsection{subsection}{2}{\z@}{-3.25ex plus-1ex
    minus-.2ex}{1.5ex plus.2ex}{\large\bf}}
\def\subsubsection{\@startsection{subsubsection}{3}{\z@}{-3.25ex plus
 -1ex minus-.2ex}{1.5ex plus.2ex}{\normalsize\bf}}

\def\eqalign#1{\,\vcenter{\openup\jot\m@th
  \ialign{\strut\hfil$\displaystyle{##}$&$\displaystyle{{}##}$\hfil
        \crcr#1\crcr}}\,}

\def\mydesc{\list{}{\labelwidth\z@ \itemindent-\leftmargin
\listparindent 1.5em
\let\makelabel\descriptionlabel}}
\let\endmydesc\endlist


\def\fnum@figure{{\small Figure \thefigure}}
\def\fakefigure{\def\@captype{figure}}

\long\def\@makecaption#1#2{
    \vskip 10pt
    \def\FCap{#2} \def\NoCap{\ignorespaces}
    \ifx \FCap\NoCap
       \setbox\@tempboxa\hbox{#1}  
      \else
       \setbox\@tempboxa\hbox{#1: \small \it #2}
    \fi
    \ifdim \wd\@tempboxa >\hsize   
        \unhbox\@tempboxa\par      
      \else                        
        \hbox to\hsize{\hfil\box\@tempboxa\hfil}
    \fi}

\pagestyle{headings}
\oddsidemargin 0.5in
\evensidemargin 0.5in
\def\@oddhead{\hbox{}\rightmark \hfil \rm\thepage}
\def\sectionmark#1{\markright {\sc{\ifnum \c@secnumdepth >\z@
      \S\thesection.\hskip 1em\relax \fi #1}}}

 \catcode`\@=12

\def\oplabel#1{
  \def\OpArg{#1} \ifx \OpArg\Empty {} \else
        \label{#1}
  \fi}

%
\newtheorem{theoremSt}{Theorem}[section]
\newtheorem{corollarySt}[theoremSt]{Corollary}
\newtheorem{propositionSt}[theoremSt]{Proposition}
\newtheorem{lemmaSt}[theoremSt]{Lemma}
\newtheorem{conjectureSt}[theoremSt]{Conjecture}
\newtheorem{exampleSt}[theoremSt]{Example}
\def\NextItem{\refstepcounter{theoremSt} \item[\thetheoremSt]}
%
\def\MakeStEnv#1{
  \newenvironment{#1}[2]{
  \begin{#1St} \oplabel{##1}%
  \global\def\CrntSt{\thetheoremSt}%
  {\rm ##2}%
}{
  \end{#1St} }
}
\MakeStEnv{theorem}
\MakeStEnv{corollary}
\MakeStEnv{proposition}
\MakeStEnv{lemma}
\MakeStEnv{conjecture}
\MakeStEnv{example}

\def\lref#1{\ref{#1} (#1)}
\newenvironment{proof}[1]{
  \def\PfArg{#1}
  \ifx\PfArg\Empty
        \edef\PfArg{\CrntSt}  \fi
 \startproof{\PfArg}%
}{
  \finishproof{\PfArg}
}
\newcommand{\startproof}[1]{
  \medbreak\mbox{}
  {\it Proof of #1:}%
}
\newcommand{\finishproof}[1]{%
\def\FPArg{#1}\ifx\FPArg\Empty\def\FPArg{\CrntSt}\fi%
\smallbreak\noindent\makebox[\textwidth]{\hfill\fbox{\FPArg}}%
\medbreak\noindent}

\newcommand{\marginwrite}[1]{}


\def\a{\alpha}
\def\eq{\!=\!}               
\def\R{{\bf R}}
\def\Z{{\bf Z}}
\def\set#1{\{ #1 \}}         



\def\today{\ifcase\month\or
   January\or February\or March\or April\or May\or June\or
   July\or August\or September\or October\or November\or December\fi
   \space\number\day, \number\year}


\title{Concentration points for Fuchsian groups}
\author{Sungbok Hong and Darryl McCullough}
\date{{\footnotesize Department of Mathematics, Korea University,
Seoul 136-701, Korea}
\\
{\footnotesize Department of Mathematics, University of Oklahoma,
Norman, OK 73019, USA}
\\
\bigskip
{\footnotesize \today}}     
\maketitle
\bigskip
{\narrower \baselineskip=14pt

\noindent{\bf Abstract:} A limit point $p$ of a discrete group of
M\"obius transformations acting on $S^n$ is called a concentration
point if for any sufficiently small connected open neighborhood $U$ of
$p$, the set of translates of $U$ contains a local basis for the
topology of $S^n$ at $p$. For the case of Fuchsian groups ($n\eq 1$),
every concentration point is a conical limit point, but even for
finitely generated groups not every conical limit point is a
concentration point. A slightly weaker concentration condition is
given which is satisfied if and only if $p$ is a conical limit point,
for finitely generated Fuchsian groups. In the infinitely generated
case, it implies that $p$ is a conical limit point, but not all
conical limit points satisfy it. Examples are given that clarify the
relations between various concentration conditions.

}

\vfil
\noindent{\footnotesize\baselineskip=12pt

\noindent AMS(MOS) Subject Classification: Primary 20H10; Secondary
30F35, 30F40, 57M50

\medskip


}

\section {Introduction}
\label{intro}

The action of a Fuchsian or Kleinian group on the sphere at infinity
can be examined from several viewpoints, and the resulting interplay
between topology, geometry, number theory, and analysis brings
richness and beauty to the subject. The topological viewpoint provides
the starting point for much of the theory, in that it gives the
dichotomy between the region of discontinuity and the limit set. The
region of discontinuity can be regarded as the portion of the sphere
at infinity with trivial or nearly trivial local dynamics. In
contrast, at points in the limit set the behavior is complicated and
varied.

For a limit point $p$ a well-known type of behavior is the property of
being a conical limit point. This property is often defined
geometrically by saying that there is a sequence of translates of the
origin (where we regard the group $\Gamma$ as acting on the Poincar\'e
ball $B^m$) that limit to $p$ and lie within a bounded hyperbolic
distance of a geodesic ray ending at $p$. But it can also be described
topologically in terms of the action of $\Gamma$ on the sphere at
infinity $S^{m-1}$. For example, one of several such characterizations
is that there exist points $q\neq r$ in $S^{m-1}$ and a sequence of
distinct elements $\gamma_n\in\Gamma$ such that $\gamma_n(p)\to q$ and
$\gamma_n(x)\to r$ for every $x\in S^{m-1}-\{p\}$. For other
topological characterizations of conical limit points,
see~\cite{A-H-M,B-M,M}.

Another topological aspect of the action of $\Gamma$ on $S^{m-1}$ is
its concentration behavior. This refers to the action of $\Gamma$ on
the set of (open) neighborhoods of $p$ in~$S^{m-1}$. The following
definitions appear in~\cite{A-H-M}.

\medskip
\noindent{\bf Definition:} An open set $U$ in $S^{m-1}$ {\it can be
concentrated} at $p$ if for every neighborhood $V$ of $p$, there
exists an element $\gamma\in\Gamma$ such that $p\in\gamma(U)$ and
$\gamma(U)\subseteq V$. If in addition the element $\gamma$ can always
be selected so that $p\in\gamma(V)$, then one says that $U$ can be
concentrated {\it with control.}
\medskip

\noindent Note that $U$ can be concentrated at $p$ if and only if the
set of translates of $U$ contains a local basis for the topology of
$S^{m-1}$ at $p$. Also, one can easily check from the definition that
(1)~there exists a neighborhood of $p$ which can be concentrated with
control if and only if there is a connected neighborhood which can be
concentrated with control (take the connected component of $U$ that
contains $p$, and require that $\gamma^{-1}(p)\in U\cap V$), and
(2)~if a neighborhood of $p$ can be concentrated with control, then
every smaller neighborhood can be concentrated with control.

\medskip
\noindent{\bf {Definition}}: The limit point $p$ is called a {\it
controlled concentration point} for $\Gamma$ if it has a neighborhood
which can be concentrated with control at~$p$.
\medskip

\noindent Concentration with control is studied in \cite{A-H-M}.
Analogously to conical limit points, $p$ is a controlled concentration
point if and only if there exist a point $r\neq p$ in $S^{m-1}$ and a
sequence $\gamma_n$ of distinct elements of $\Gamma$ so that
$\gamma_n(p)\to p$ and $\gamma_n(x)\to r$ for all $x\in
S^{m-1}-\{p\}$. In particular, every controlled concentration point is
a conical limit point. However, examples are given in \cite{A-H-M} of
conical limit points of 2-generator Schottky groups which are not
controlled concentration points (see also proposition~\ref{4.1}
below). For groups of divergence type, controlled concentration points
have full Patterson-Sullivan measure in the limit set. There is a
direct connection between controlled concentration points and the
dynamics of geodesics in the hyperbolic manifold $B^m/\Gamma$. Call a
geodesic ray in $B^m/\Gamma$ {\it recurrent} if it is the image of a
geodesic ray in $B^m$ that ends at a controlled concentration point.
In an appropriate metric, the space of recurrent geodesic rays in
$B^m/\Gamma$ is a metric completion of the space of closed geodesics
in $B^m/\Gamma$ (where both spaces are topologized as subspaces of the
unit tangent bundle of~$B^m/\Gamma$).

We turn now to weaker concentration properties. It is not difficult to
show (see \cite{M1}) that every limit point $p$ has a disconnected
neighborhood that can be concentrated at $p$. So the weakest
reasonable concept of concentration behavior is the following.

\medskip
\noindent{\bf {Definition}}: The limit point $p$ is called a  {\it
weak concentration point} for $\Gamma$ if there exists a connected
open set that can be concentrated at~$p$.
\medskip

\noindent Weak concentration points are studied in \cite{M1}. It turns
out that for a geometrically finite group, every limit point is a weak
concentration point, and for any group, all but countably many limit
points are weak concentration points. A more restrictive condition is
that {\it every} sufficiently small connected neighborhood can be
concentrated:

\medskip
\noindent{\bf {Definition}}: The limit point $p$ is called a  {\it
concentration point} for $\Gamma$ if every sufficiently small
connected neighborhood of $p$ can be concentrated at~$p$.
\medskip

\noindent In this paper, we will study concentration behavior for
Fuchsian groups. From now on, let $\Gamma$ be Fuchsian. A slightly
weaker concept than concentration point turns out to be important.

\medskip
\noindent{\bf Definition:} The limit point $p$ is called a {\it
geodesic separation point} for the Fuchsian group $\Gamma$ if for
every sufficiently small connected neighborhood $U$ of $p$, either $U$
or $S^1-\overline{U}$ can be concentrated at~$p$.
\medskip

\noindent The name of this property derives from the fact that for a
geodesic separation point $p$, if $\lambda$ is any geodesic in $B^2$
whose endpoints separate $p$ from the boundary of a small neighborhood
of $p$, then for any neighborhood $V$ of $p$ there exists
$\gamma\in\Gamma$ so that the endpoints of $\gamma(\lambda)$ separate
$p$ from the boundary of~$V$. Indeed, it is easily verified that this
is equivalent to the condition in the definition; this simply uses the
fact that every connected neighborhood of $p$ (other than $S^1$
itself) is an interval, so corresponds to the unique geodesic in $B^2$
that runs between its endpoints.

The purpose of this paper is to investigate the general relations
between these concentration properties for Fuchsian groups. We
summarize them here; unless otherwise stated, references are to
results that appear later in this paper. By $\Gamma$ we denote an
nonelementary Fuchsian group, by $\Gamma_0$ a certain two-generator
Schottky group defined in~\S4, and by $\Gamma_1$ a certain infinitely
generated Fuchsian group defined in~\S5.

\begin{enumerate}
\item[{\rm1.}] Every controlled concentration point for $\Gamma$ is a
concentration point (a direct consequence of the definitions). There
are uncountably many concentration points for $\Gamma_0$ that are not
controlled concentration points (theorem~\ref{4.2}).
\item[{\rm2.}] Every concentration point for $\Gamma$ is a geodesic
separation point (immediate from the definitions). There are
uncountably many geodesic separation points for $\Gamma_0$ that are
not concentration points (theorem~\ref{4.3} and theorem~\ref{3.1}).
\item[{\rm3.}] Every geodesic separation point for $\Gamma$ is a
conical limit point (proposition~\ref{5.2}). There are uncountably
many conical limit points for $\Gamma_1$ which are not geodesic
separation points (proposition~\ref{5.1}). However, if $\Gamma$ is
finitely generated, then every conical limit point is a geodesic
separation point (theorem~\ref{3.1}).
\item[{\rm4.}] Every conical limit point or parabolic fixed point for
$\Gamma$ is a weak concentration point (theorem~3.1 of~\cite{M1}).
There are uncountably many weak concentration points for $\Gamma_1$
which are neither conical limit points nor parabolic fixed points
(proposition~\ref{5.1}). However, if $\Gamma$ is finitely generated,
then every weak concentration point is either a conical limit point or
a parabolic fixed point (by theorem~2 of~\cite{B-M}).
\item[{\rm5.}] At most countably many limit points of $\Gamma$ are not
weak concentration points (theorem~3.6 of~\cite{M1}). However, if
$\Gamma$ is finitely generated, then every limit point is a weak
concentration point (corollary~2.2 of~\cite{M1}).
\end{enumerate}

We assume familiarity with the basic concepts of M\"obius groups as
exposited, for example, in \cite{B}. We use the term {\it Nielsen
hull} for the (hyperbolic) convex hull in $B^2$ of the limit set of a
Fuchsian group $\Gamma$, and {\it Nielsen core} for the quotient of
the Nielsen hull by $\Gamma$. Otherwise, our terminology and notation
are standard. The reader may find it useful to examine the examples of
\S\S4 and 5 before delving into \S\S2 and~3, whose main objective is
the proof of theorem~\ref{3.1}.

This manuscript is a revised version of a preprint that was circulated
in 1992. The proof of theorem~\ref{3.1} is substantially rewritten,
and a gap in it has been filled. Section~5 is new, and the
introduction has been completely rewritten to place the results in the
context of subsequent developments which have appeared
in~\cite{A-H-M,M1}.

\section [Controlled Concentration Points] {Controlled  concentration
points and geodesic laminations}
\label{ccp}
\medskip

This section contains some sufficient conditions for a limit point of
a Fuchsian group to be a controlled concentration point. In
particular, if $L$ is a compact geodesic lamination in a hyperbolic
$2$-manifold, then every endpoint of a leaf of the preimage of $L$ in
$B^2$ is a controlled concentration point.

Our first condition follows from theorem~2.3 of \cite{A-H-M}, but for
simplicity we give a direct self-contained argument.

\begin{lemma}{2.1}{}
Let $\Gamma$ be a torsionfree discrete group of M\"obius
transformations acting on the Poincar\'e disc $B^m$, and let
$\pi\colon B^m\to B^m/\Gamma$ be the quotient map. Let $y_0\in B^m$
and let $p$ be a point in $S^{m-1}$. Let $\a\colon[0,\infty)\to B^m$
be the geodesic ray from $y_0$ to $p$, parameterized at unit speed.
Suppose further that there exist numbers $t_i$, with $\a(t_i)$
limiting to $p$, so that in the tangent bundle $T(B^m/\Gamma)$, the
images $d\pi(\a'(t_i))$ converge to $d\pi(\a'(0))$. Then $p$ is a
controlled concentration point for~$\Gamma$.
\marginwrite{2.1}
\end{lemma}

\begin{proof}{} Consider the hyperbolic codimension~1 hyperplane
through $y_0$ perpendicular to $\alpha$, and let $U$ be the
neighborhood of $p$ in $S^{m-1}$ which is one of the components
of the complement of the boundary of the hyperplane.  For fixed
positive $n$, let $V_n$ be the smaller neighborhood of $p$ whose
boundary is the boundary of the hyperplane perpendicular to $\alpha$
and crossing it at hyperbolic distance $n$ from $y_0$. The hypothesis
of lemma~2.1 shows there are elements $\gamma_i\in\Gamma$ which
translate $\alpha'(0)$ close to $\alpha'(t_i)$. For all sufficiently
large $i$, the geodesic $\gamma_i\circ\alpha$ follows along very close
to $\alpha$ for more than distance $n$ --- close enough so that $p\in
\gamma_i(V_n)$. By making $i$ perhaps even larger, one can also ensure
that $\gamma_i(U)\subseteq V_n$, showing that $p$ is a controlled
concentration point.
\end{proof}

For a Fuchsian group $\Gamma$ acting on the Poincar\'e disc $B^2$, we
denote the Nielsen core of $B^2/\Gamma$ by $N(B^2/\Gamma)$. If
$\Gamma$ is not elementary, then the interior of $N(B^2/\Gamma)$ is
nonempty. Since the Nielsen hull is convex, it follows that if a
geodesic in $B^2/\Gamma$ leaves the Nielsen core, it will never
reenter. Note that such a geodesic cannot lie in a compact subset of
$B^2/\Gamma$. When $\Gamma$ is finitely generated and torsionfree, its
Nielsen core is a 2-manifold of finite area whose boundary is a finite
collection of simple closed geodesics.

\begin{theorem}{2.2}{}
Let $\Gamma$ be a Fuchsian group acting on the Poincar\'e disc $B^2$.
Suppose there exists a geodesic ray in $B^2$ which ends at $p\in S^1$,
which has no transverse crossing with any of its translates, and whose
image in $B^2/\Gamma$ lies in a compact subset. Then $p$ is a
controlled concentration point for~$\Gamma$.
\marginwrite{2.2}
\end{theorem}

\begin{proof}{} Let $S$ denote $B^2/\Gamma$, let $\alpha$
denote the hypothesized geodesic ray, and let $\alpha_0$ denote its
image in $S$. Now $S$ is an increasing union of compact suborbifolds,
and by filling in any complementary discal 2-orbifolds for these
suborbifolds, we may assume that each has orbifold fundamental group
which injects into $\Gamma$. Replacing $\Gamma$ by the fundamental
group of one of the suborbifolds that contains $\alpha_0$, we may assume
that $\Gamma$ is finitely generated. Passing to a subgroup of finite
index, we may assume that $\Gamma$ is torsionfree.

We will use the elementary theory of geodesic laminations as presented
in \cite{C-B}. Let $L$ be the set of points $y\in S$ with the
following property: there is a sequence of points $x_i$ on
$\alpha$ that limits to $p$, whose images limit to $y$.
Since $\alpha_0$ has no transverse self-intersections, it follows that
$L$ is a nonempty compact geodesic lamination in $S$, for which each
tangent vector is a limit of a sequence of tangent vectors to $\alpha_0$
at points whose preimages limit to~$p$. Since $L$ is compact, it must
lie in~$N(S)$.

Suppose that $L$ contains a simple closed geodesic $C$. If there were
no collar neighborhood of $C$ that contained the image of a subray of
$\alpha$, then $\alpha_0$ would have transverse self-intersections.
Therefore $\alpha_0$ must either be contained in $C$, or spiral toward
$C$. In either case, $p$ is the endpoint of the axis of a hyperbolic
translation in $\Gamma$, so is a controlled concentration point.

Suppose that $L$ does not contain any simple closed geodesics, and
that $\alpha_0$ is contained in a leaf of $L$. Orient $\alpha$ in the
direction pointing toward $p$, and let $v$ be the initial (unit)
tangent vector of $\alpha_0$. Let $x_i$ be a sequence of points of
$\alpha$, limiting to $p$, whose images $y_i$ limit to the starting
point of $\alpha_0$, and let $v_i$ be the oriented tangent vectors to
$\alpha_0$ at $y_i$. If a subsequence of the  $v_i$ limits to $v$,
then lemma~2.1 shows that $p$ is a controlled concentration point. So
we may assume that they limit to $-v$. Let $\beta$ be the subarc of
$\alpha_0$ from its initial point to some $y_i$. For $j$ sufficiently
large so that $v_j$ is extremely close to $-v$, the part of $\alpha_0$
that ends at $y_j$ follows backwards along $\beta$ staying very close
until after it passes $y_i$, and at a point where it passes $y_i$ it
is pointing approximately in the direction of $v$ (see Figure~1).
Therefore lemma~2.1 still applies.

\begin{figure}[htb]

\centerline{
\epsfig{figure=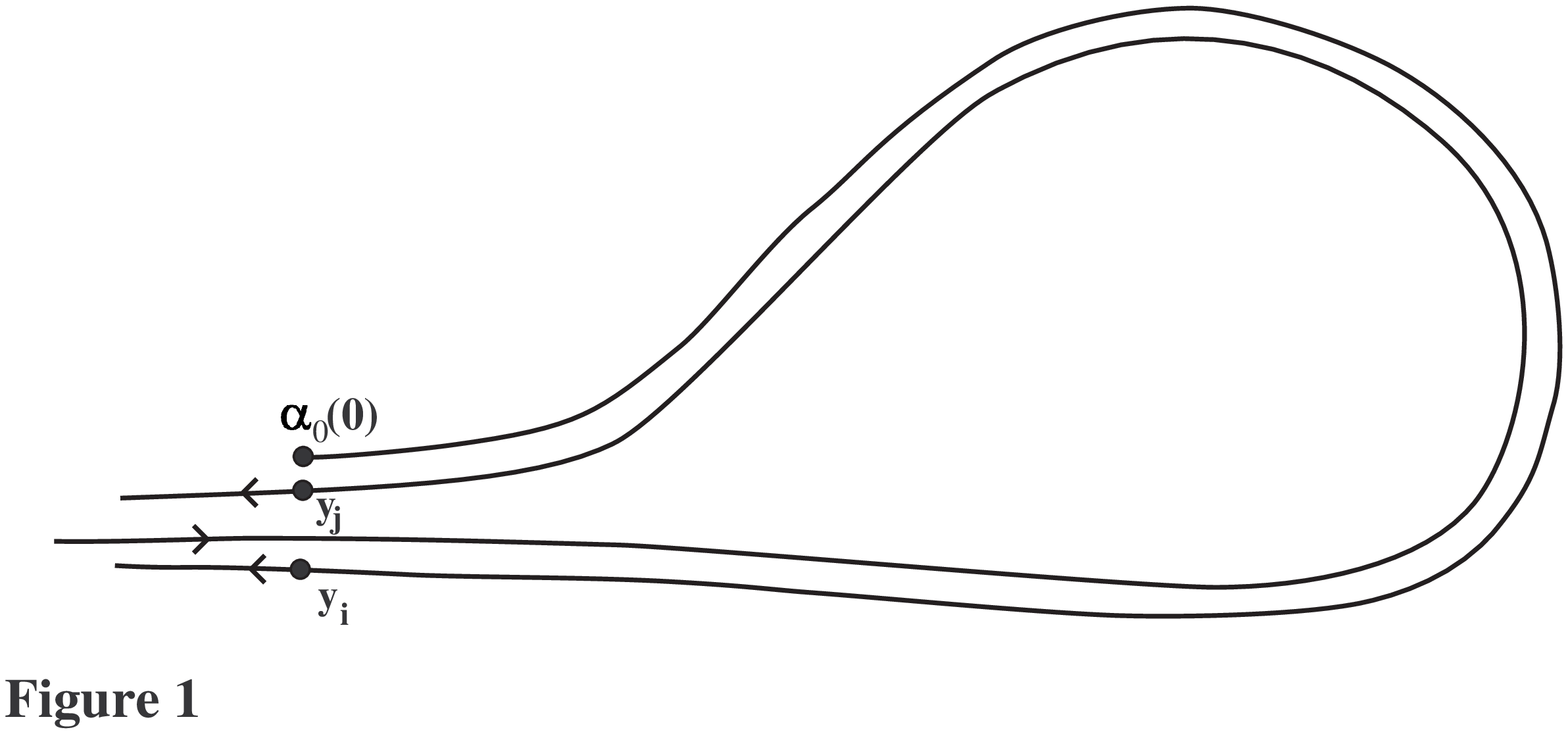,width=6cm,height=10cm,angle=-90}}
\end{figure}

The last case is that $L$ does not contain any simple closed
geodesics, and $\alpha_0$ does not lie in a leaf of $L$. Note that
$\alpha_0$ must be disjoint from $L$, since any transverse crossings
with $L$ would force self-intersections of $\alpha_0$. A {\it
principal region} for $L$ is a component of $N(S)-L$. Since $L$ lies
in $N(S)$ and does not contain any of the boundary geodesics, we may
assume (by shortening $\alpha$) that $\alpha_0$ lies in $N(S)$. Let
$U$ be the principal region for $L$ that contains $\alpha_0$.
Following lemma~4.4 of \cite{C-B}, we wish to describe the
possibilities for $U$. In that lemma, $S$ is closed and has no cusps,
and $U$ is either isometric to the interior of a finite-sided ideal
polygon, or there is a compact submanifold $U_0$ of $U$, whose
boundary is a union of closed geodesics, such that $U-U_0$ is
isometric to the interior of a finite collection of crowns. (A {\it
crown} is a complete hyperbolic surface with finite area and geodesic
boundary, which is homeomorphic to $S^1\times[0,1]-Q$ for some finite
subset $Q$ of $S^1\times\set{1}$.) We define a {\it cuspidal crown} to
be a complete hyperbolic surface with finite area and geodesic
boundary, which is homeomorphic to $(B^2\cup \partial
B^2)-(\set{0}\cup Q)$, where $Q$ is a finite subset of~$\partial B^2$.
For later reference, we state the next observation as a lemma.

\begin{lemma}{Casson-Bleiler}{} Let $L$ be a geodesic lamination in a
hyperbolic 2-manifold $F$ of finite area, with boundary consisting of
closed geodesics, and let $U$ be a component of $F - L$. Then either
\begin{enumerate}
\item[{\rm1.}]
$U$ is isometric to the interior of a finite-sided ideal polygon in
$B^2$, or
\item[{\rm2.}]
there is a submanifold $U_0$ of $U$, whose boundary consists of closed
geodesics, such that $U-U_0$ is isometric to the interior of a finite
collection of crowns, or
\item[{\rm3.}]
$U$ is isometric to the interior of a cuspidal crown.
\end{enumerate}
\marginwrite{Casson-Bleiler}
\end{lemma}

\begin{proof}{}
The proof is exactly like the proof of lemma~4.4 of \cite{C-B}, with
the case (3) arising when the element called $g$ there is parabolic.
\end{proof}

Returning to the proof of theorem~\ref{2.2}, since $\alpha_0$ is
disjoint from $L$, but limits onto $L$, it must limit onto one of the
noncompact boundary geodesics of a principal region of one of the
forms described in lemma~\ref{Casson-Bleiler}. It follows that there
is a leaf in the preimage of the boundary leaves of $L$ that ends at
$p$. Replacing $\alpha$ by a subray of that leaf ending at $p$, we are
in the previous case where $\alpha_0$ was assumed to lie in a leaf,
and again it follows that $p$ is a controlled concentration point.
\end{proof}

\begin{corollary}{2.3}{}
Let $\Gamma$ be a torsionfree Fuchsian group, and let $L$ be a compact
geodesic lamination in $B^2/\Gamma$. Then the endpoints of the leaves
of the preimage of $L$ in $B^2$ are controlled concentration points
for~$\Gamma$.
\marginwrite{2.3}
\end{corollary}

\begin{proof}{} Apply theorem 2.2 to a subray of the leaf.
\end{proof}

\section [Geodesic Separation Points] {Geodesic separation points}
\label{gsp}

The first result of this section, proposition~\ref{5.2}, shows that
every geodesic separation point is a conical limit point. In
particular, a parabolic fixed point cannot be a geodesic separation
point. On the other hand, in \S5 we will give an example of an
infinitely generated Fuchsian group with uncountably conical limit
points which are not geodesic separation points. Theorem~\ref{3.1}
shows that this cannot happen in a finitely generated example, since
then every limit point is either a parabolic fixed point or a geodesic
separation point.

Note that proposition~\ref{5.2} implies that for Fuchsian groups,
concentration points are conical limit points. Whether this holds in
higher dimensions is an open question.

\begin{proposition}{5.2}{}
Let $p$ be a geodesic separation point of a Fuchsian group $\Gamma$.
Then $p$ is a conical limit point.
\marginwrite{5.2}
\end{proposition}

\begin{proof}{}  Clearly, $p$ cannot be the endpoint of an interval of
discontinuity with finite stabilizer. If $p$ is the endpoint of an
interval of discontinuity with infinite stabilizer, then $p$ is the
attracting fixed point of a hyperbolic element and hence is a conical
limit point. So we may assume that every neighborhood of $p$ contains
limit points of $\Gamma$ on both sides of $p$. Let $W$ be a
neighborhood of $p$ such that for every connected neighborhood $U$ of
$p$ with $U\subseteq W$, either $U$ or $S^1-\overline{U}$ can be
concentrated at $p$. By the Double Density Theorem, there exists an
axis $\lambda$ of a hyperbolic element of $\Gamma$ whose endpoints lie
in $W$ and separate $p$ from the boundary of $W$. If $x$ is a point on
this axis, then translates of $x$ lie at intervals of some fixed
length $d$ along $\lambda$. Since $p$ is a geodesic separation point,
there must be translates of $\lambda$ that intersect $\a$ arbitrarily
close to $p$. On each such translate of $\lambda$, there are
translates of $x$ within hyperbolic distance $d$ of $\a$. Therefore
$p$ is a conical limit point.
\end{proof}

Several times in the proof of theorem~\ref{3.1}, we will use the
observation that if a portion of a geodesic ray in $N(B^2/\Gamma)$
moves far out a cusp, but the ray does not continue all the way out
the cusp, then it must behave as follows. For some time it travels
almost straight out the cusp, then it starts to spiral around the
cusp, finally becoming tangent to some horocycle, then it returns to
the thick part of $N(B^2/\Gamma)$. (This behavior is easily seen by
normalizing so that the parabolic element generating the cusp acts in
the upper half-plane model as $z\mapsto z+1$, and observing the
behavior of a geodesic arc that rises to a high vertical coordinate
before descending to the real line.) Note that in particular, any such
ray in $N(B^2/\Gamma)$ must have self-intersections, and when it is
spiraling near its tangent horocycle it must intersect any geodesic
ray that travels a great deal further out the cusp at nearly right
angles.

\begin{theorem}{3.1}{}
Let $\Gamma$ be a Fuchsian group. If $\Gamma$ is finitely generated,
then every limit point of $\Gamma$ is either a parabolic fixed point
or a geodesic separation point.
\marginwrite{3.1}
\end{theorem}

\begin{proof}{} Assume that $p$ is not a parabolic fixed point.
Passing to a subgroup of finite index, we may assume that $\Gamma$ is
torsionfree. If $\Gamma$ is elementary, then $p$ is the endpoint of
the axis of a hyperbolic element, and hence a controlled concentration
point. So we will assume that $\Gamma$ is nonelementary, and hence
that its Nielsen core has nonempty interior. Since $\Gamma$ is
finitely generated, its Nielsen core has finite area and has boundary
a (possibly empty) finite collection of simple closed geodesics.

For later reference, we isolate the next argument as a lemma.

\begin{lemma}{3.2}{}
Let $\Gamma$ be a finitely generated torsionfree Fuchsian group, and
let $p$ be a limit point of $\Gamma$ which is not a controlled
concentration point and is not a parabolic fixed point. Let $\alpha$
be a geodesic ray ending at $p$, and lying in the interior of the
Nielsen hull of $\Gamma$. Suppose $\lambda$ is a geodesic in $B^2$
that intersects $\alpha$, such that $p$ is not a limit point of the
crossings of the translates of $\lambda$ with $\alpha$. Then there
exists a finitely generated subgroup $\Gamma'$ of $\Gamma$ with the
following properties:
\begin{enumerate}
\item[{\rm(a)}] $p$ is a limit point of $\Gamma'$, and some subray of
$\alpha$ lies in the interior of the Nielsen hull of $\Gamma'$.
\item[{\rm(b)}] Either the area of $N(B^2/\Gamma')$ is less than the
area of $N(B^2/\Gamma)$, or the areas are equal and the genus of
$N(B^2/\Gamma')$ is less than the genus of~$N(B^2/\Gamma)$.
\end{enumerate}
\marginwrite{3.2}
\end{lemma}

\begin{proof}{} Since $\alpha$ lies in the Nielsen hull of $\Gamma$,
its image $\alpha_0$ in $B^2/\Gamma$ lies in $N(B^2/\Gamma)$. Let
$\Gamma\lambda$ be the union of the translates of $\lambda$. Then some
subray of $\alpha$ lies entirely in a component $E$ of
$B^2-\Gamma\lambda$. Note that $E$ is convex, since it is an
intersection of half-spaces. Let $\Gamma'$ be the stabilizer of $E$.
Since $E$ is precisely invariant under $\Gamma$, $E/\Gamma'$ maps
injectively into $B^2/\Gamma$ under the covering projection
$B^2/\Gamma'\rightarrow B^2/\Gamma$.

Suppose for contradiction that $p$ is not a limit point of $\Gamma'$.
Then some subray of $\alpha$ injects into $E/\Gamma'$, and hence into
$N(B^2/\Gamma)$. If $\alpha_0$ does not lie in a compact subset of the
Nielsen core, then since $p$ is not a parabolic fixed point,
$\alpha_0$ must enter and leave some cusp of $N(B^2/\Gamma)$
infinitely many times, moving farther and farther out toward the end
of the cusp. This is impossible as a subray of $\a$ maps injectively.
So $\alpha_0$ lies in a compact subset. But then, theorem~2.2 implies
that $p$ is a controlled concentration point for $\Gamma$, giving the
contradiction.

Now $\Gamma'$ is nonelementary, since otherwise $p$ would be a fixed
point of a hyperbolic element, and hence a controlled concentration
point. Therefore the interior of the Nielsen hull of $\Gamma'$ is
nonempty. Since the orbit of any point of $E$ under $\Gamma'$ lies in
$E$, the limit set of $\Gamma'$ lies in $\overline{E}\cap S^1$. Since
$E$ is convex, this shows that the interior of the Nielsen hull of
$\Gamma'$ lies in $E$. Therefore the covering map from $B^2/\Gamma'$
to $B^2/\Gamma$ restricts to a map $N(B^2/\Gamma')\to N(B^2/\Gamma)$
whose restriction to the interior of $N(B^2/\Gamma')$ is an isometric
imbedding. In particular, this shows that $\Gamma'$ is finitely
generated, so that $N(B^2/\Gamma')$ has boundary consisting of closed
geodesics. Consequently, the topological effect of the map
$N(B^2/\Gamma')\to N(B^2/\Gamma)$ must be first to include
$N(B^2/\Gamma')$ into a larger (or possibly equal) surface, then
(possibly) identify some pairs of boundary components. Therefore
either $N(B^2/\Gamma')$ has smaller area than $N(B^2/\Gamma)$, or it
has the same area and has smaller genus, or the restriction of the
covering map is a homeomorphism from $N(B^2/\Gamma')$ to
$N(B^2/\Gamma)$. In the latter case, $\Gamma\eq\Gamma'$ so their
Nielsen hulls are equal. But the interior of the Nielsen hull of
$\Gamma'$ is disjoint from the translates of $\lambda$, so $\a$ could
not have intersected $\lambda$. Therefore $\Gamma'\eq\Gamma$ is
impossible, giving assertion~(b).

Finally, if no subray of $\alpha$ lies in the interior of the convex
hull of $\Gamma'$, then since $\alpha$ ends at a limit point of
$\Gamma'$, $\alpha_0$ must either coincide with or spiral onto a
boundary geodesic for $N(B^2/\Gamma')$. But then $p$ is a fixed point
of a hyperbolic element of $\Gamma'$, a contradiction. This completes
the proof of assertion~(a).
\end{proof}

We now continue the proof of theorem~\ref{3.1}. Suppose for
contradiction that $p$ is not a geodesic separation point. Fix a
geodesic ray $\alpha$ in $B^2$ running from a point in the interior of
the Nielsen hull of $\Gamma$ to $p$, and let $\alpha_0$ denote its
image in $N(B^2/\Gamma)$. Make an initial choice of connected
neighborhood $W$ of $p$, small enough so that whenever the endpoints
of $\lambda$ lie in $W-\set{p}$, they separate $p$ from the boundary
of $W$ if and only if $\lambda$ intersects $\alpha$. Since $p$ is not
a geodesic separation point, there exists a geodesic $\lambda$ with
endpoints in $W-\set{p}$, so that $\lambda$ intersects $\alpha$ but
for which there is a neighborhood $V$ of $p$ for which no translate of
$\lambda$ separates $p$ from the boundary of~$V$.

Suppose first that $p$ is not a limit point of the intersections of
the translates of $\lambda$ with $\alpha$. Let $\Gamma'$ be a subgroup
of $\Gamma$ obtained using lemma~\ref{3.2}. By condition~(a) of
lemma~\ref{3.2}, $p$ is a limit point of $\Gamma'$, and since
$\Gamma'$ is a subgroup of $\Gamma$, $p$ is not a geodesic separation
point for $\Gamma'$. Replace $\Gamma$ by $\Gamma'$, shorten $\alpha$
if necessary so that it lies in the interior of the Nielsen hull of
$\Gamma'$, and replace $W$ by a smaller neighborhood if necessary. By
condition~(b) of lemma~\ref{3.2}, such a procedure can only occur
finitely many times. So we may assume that $p$ is a limit point of the
intersections of the translates of $\lambda$ with~$\alpha$.

Let $x_i$ be a sequence of intersection points of translates
$\gamma_i(\lambda)$ with $\alpha$, which limit to $p$. Since $p$ is
not a geodesic separation point, the angles between $\alpha$ and
$\gamma_i(\lambda)$ at the $x_i$ must limit to $0$. Moreover, by
passing to a subsequence we may assume that for one of the endpoints
$e_1$ of $\lambda$, the sequence $\gamma_i(e_1)$ limits to $p$, while
for the other endpoint $e_2$, the sequence $\gamma_i(e_2)$ limits to a
point $q$ distinct from~$p$ (since the $\gamma_i(e_2)$ lie in
$S^1-\overline{V}$).

Suppose for contradiction that no subsequence of the images of the
$x_i$ in $N(B^2/\Gamma)$ converges. By passing to a subsequence, we
may assume that for some cusp of $N(B^2/\Gamma)$, the images of the
$x_i$ lie farther and farther out the cusp. Now $\alpha_0$ does not
travel straight out the cusp, since $p$ is not a parabolic fixed
point. Therefore there are portions of $\alpha_0$ that travel almost
straight out the cusp for a long time, then start to spiral, becoming
tangent to some horocycle, then return back to the thick part of
$B^2/\Gamma$. The image $\lambda_0$ of $\lambda$ either travels
straight out the cusp, or has infinitely many portions similar to
those of $\alpha_0$. Consider such a portion of $\alpha_0$. At the
part where it spirals near its tangent horocycle, it crosses
$\lambda_0$ almost perpendicularly (either where $\lambda_0$ is
traveling straight out the cusp, or on infinitely many portions that
are traveling out to horocycles much farther out the cusp). This
implies that there is a sequence of nearly perpendicular crossings of
$\alpha$, converging to $p$. This contradicts the choice of $\lambda$.
Therefore, by taking a subsequence of the $x_i$, we may assume that
the images of the $x_i$ converge to a point~$s$ in $N(B^2/\Gamma)$,
and moreover that the images of the (unit) tangent vectors to
$\alpha$ at the $x_i$ (oriented to point toward $p$) also converge.
Since the intersection angle between $\alpha$ and the
$\gamma_i(\lambda)$ at the $x_i$ limit to $0$, the images of the
tangent vectors of the $\gamma_i(\lambda)$ at the $x_i$ (oriented to
point toward the endpoint $\gamma_i(e_1)$ of $\gamma_i(\lambda)$) must
also converge to the same limiting vector. Let $\mu_0$ be the geodesic
in $B^2/\Gamma$ determined by this tangent vector.

Suppose for contradiction that $\mu_0$ has a transverse
self-intersection. Then some segment $\sigma$ of $\mu_0$ containing
$s$ has a transverse self-intersection at a point $s_0$, where it
crosses itself making a positive intersection angle $\theta$. There
are portions of $\alpha_0$ and $\lambda_0$ which approximate $\sigma$
arbitrarily closely. Therefore there are crossings of $\alpha_0$ with
$\lambda_0$ close to $s_0$, at angles close to $\theta$. This shows
that there are translates of $\lambda$ crossing $\alpha$ at angles
approximately $\theta$, with the crossings limiting to $p$. This
contradicts the choice of $\lambda$. We conclude that $\mu_0$ has no
transverse self-intersections. By similar reasoning, $\lambda_0$
cannot intersect~$\mu_0$ transversely.

Every point of $\mu_0$ must be a limit of points of $\alpha_0$, and
hence $\mu_0$ lies in~$N(B^2/\Gamma)$. Suppose for contradiction that
it does not lie in a compact subset of $N(B^2/\Gamma)$. Since it has
no self-intersections, it must travel all the way out a cusp of
$N(B^2/\Gamma)$. Since $\lambda_0$ does not cross $\mu_0$
transversely, but every point of $\mu_0$ is a limit of points of
$\lambda_0$, $\lambda_0$ must also travel straight out that cusp.
Since $\alpha_0$ cannot travel straight out the cusp, because $p$ is
not a parabolic fixed point, there must as before be infinitely many
portions of $\alpha_0$ that spiral near horocycles in this cusp. This
produces nearly perpendicular intersections of $\alpha_0$ with
$\lambda_0$, as before contradicting the choice of $\lambda$. So we
may assume that $\mu_0$ lies in a compact subset of $N(B^2/\Gamma)$.
Since $\mu_0$ has no self-intersections, its closure is a geodesic
lamination $L$, for which every tangent vector is a limit of tangent
vectors of~$\alpha_0$.

Suppose for contradiction that every subray of $\alpha_0$ has
transverse crossings with $\mu_0$. Let $\mu$ be a lift of $\mu_0$ to
$B^2$. Then there is a sequence of translates of $\mu$ intersecting
$\alpha$ in a sequence of points $r_i$ converging to $p$. By passing
to subsequences, we may assume that $r_i$ and $x_i$ alternate as one
moves along $\alpha$. Since $\lambda_0$ does not cross $\mu_0$
transversely, the translates must be disjoint and must alternate as
shown in Figure~2. Therefore there is a sequence of translates of
$\mu$ that converges to the geodesic from $p$ to $q$. Since $L$ is
closed, this shows that the geodesic from $p$ to $q$ is the lift of a
geodesic of $L$. By corollary~\ref{2.3}, $p$ is a controlled
concentration point for $\Gamma$, a contradiction. So by passing to a
subray of $\alpha$, we may assume that $\alpha_0$ does not cross $L$
transversely. Since $p$ is not a controlled concentration point,
corollary~\ref{2.3} shows that $\alpha_0$ does not lie in a leaf of
$L$, so $\alpha_0$ is disjoint from~$L$.

\begin{figure}[htb]
\centerline{
\epsfig{figure=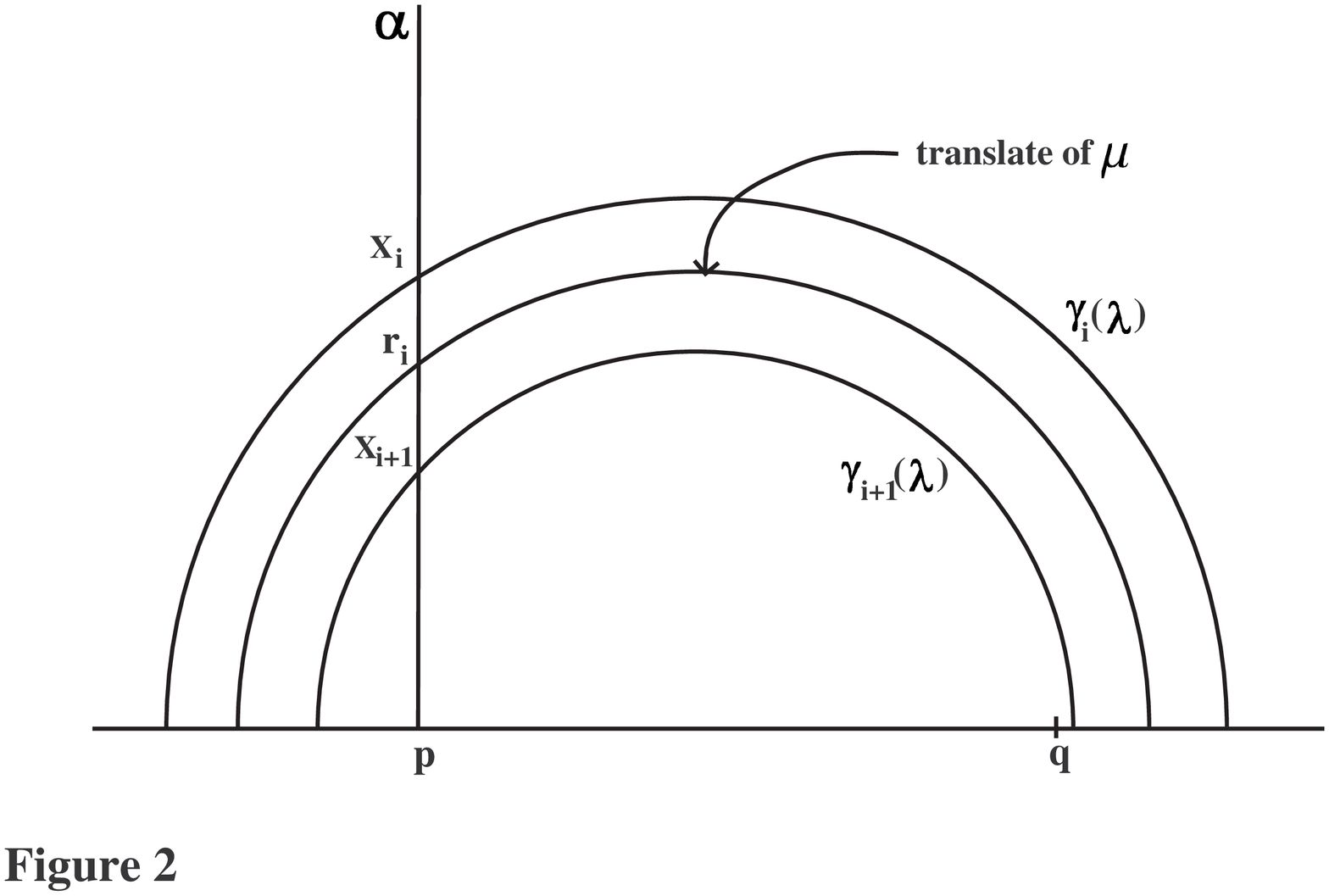,width=8cm,height=9cm,angle=-90}}
\end{figure}

Cutting $N(B^2/\Gamma)$ along $L$, we obtain pieces as described in
lemma~\ref{Casson-Bleiler}. Now $\alpha_0$ lies in one of these pieces
and tangent vectors of $\alpha_0$ lie arbitrarily close to vectors in
a boundary geodesic $\rho_0$ of $L$. Suppose this geodesic is not
closed. Then it lies in a polygon or a crown, and some terminal
segment of $\alpha_0$ travels out an end of the polygon or crown
limiting onto $\rho_0$. It follows that some lift of $\rho_0$ ends at
$p$. By corollary~\ref{2.3}, $p$ is a controlled concentration point,
a contradiction. So we may assume that $\rho_0$ is a simple closed
loop.

We now refer to Figure~3. Since tangent vectors of $\lambda_0$ limit
to $\rho_0$, there is a portion of $\lambda_0$ that spirals very close
to $\rho_0$, reaches a minimum distance $d$, then spirals away.
Assuming that the portion is selected to make $d$ sufficiently small,
there is a sequence of lifts of $\lambda_0$ to $B^2$ that appear with
a lift $\rho$ of $\rho_0$ as shown. Similarly, there are portions of
$\alpha_0$ that spiral in toward $\rho_0$, reaching a minimum distance
$d'$ from $\rho_0$, where $d'$ may be selected to be much smaller than
$d$, and then spiral away. This implies there are translates $\alpha$
as shown in Figure~3. Notice that any such lift must cross one of the
translates of $\lambda$ in Figure~3, making an angle larger than some
positive lower bound $\theta_0$. A succession of translates of
$\alpha$ corrsponding to smaller and smaller values of $d'$ shows that
there is a sequence of intersections of $\alpha$ with translates of
$\lambda$, converging to $p$, at which the crossing angles are all
greater than $\theta_0$. This again contradicts the choice of
$\lambda$, establishing that $p$ is a geodesic separation point.
\end{proof}

\begin{figure}[htb]

\centerline{
\epsfig{figure=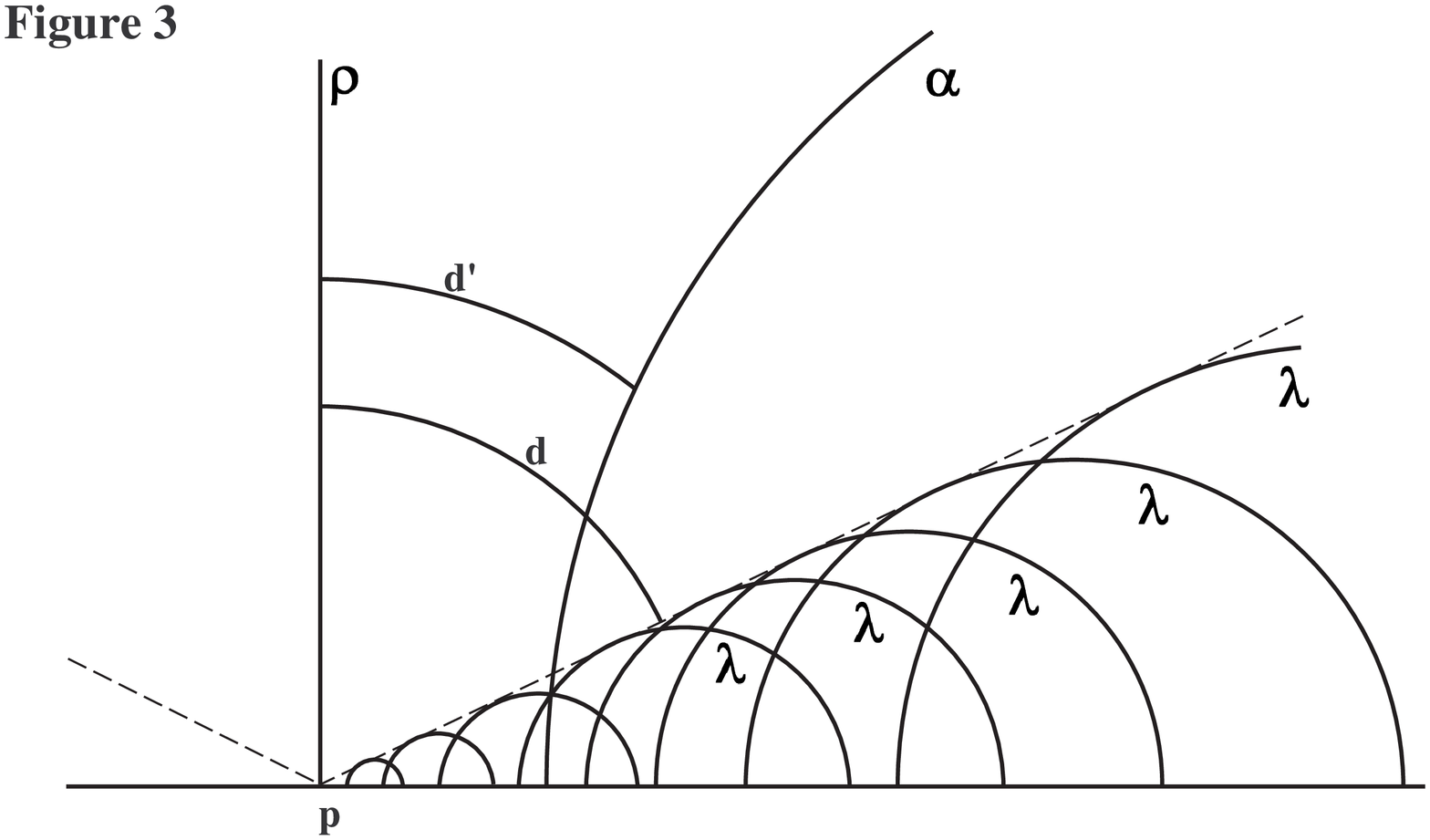,width=9cm,height=11.5cm,angle=-90}}
\end{figure}

\section [The Schottky Example]{The Schottky example}
\label{examples}

In this section, we show that the sets of controlled concentration
points, concentration points, and geodesic separation points can
differ even for finitely generated Fuchsian groups. For simplicity, we
will work with an explicit two-generator 2-dimensional Schottky group
$\Gamma_0$, although it will be apparent that the same phenomena occur
for other similar examples. The limit set of $\Gamma_0$ is a Cantor
set which can be understood quite explicitly using the sequence of
crossings of a geodesic ray (ending at the limit point) with the
translates of two fixed geodesics which lie in the boundary of a
fundamental domain.

To define $\Gamma_0$, we work in the Poincar\'e unit disc $B^2$.
Figure~4\ shows a fundamental domain for the action of $\Gamma_0$ on
$B^2$. Its frontier has two geodesics, $a$ and $a'$, with centers on
the real axis, and two more, $b$ and $b'$, with centers on the
imaginary axis. It is generated by two isometries, one of which
preserves the real axis and carries $a'$ to $a$, and the other
preserving the imaginary axis and carrying $b'$ to $b$. Fix
arbitrarily a normal direction to $a$ to call the positive direction.
It determines a normal direction for each translate of $a$; the side
of the translate into which it points will be called the {\it positive
side.} Similarly, we select a positive side for $b$ and its
translates. A crossing of an oriented geodesic or geodesic ray in
$B^2$ with a translate of $a$ or $b$ will be called a {\it positive}
crossing when it crosses from the negative side to the positive side,
otherwise it will be called a {\it negative} crossing.

A few more translates of these geodesics are drawn in Figure~4; it is
convenient to label all translates of $a$ with a letter $a$, located
{\it on the positive side,} and similarly to label all translates
of~$b$.
\begin{figure}[htb]

\centerline{
\epsfig{figure=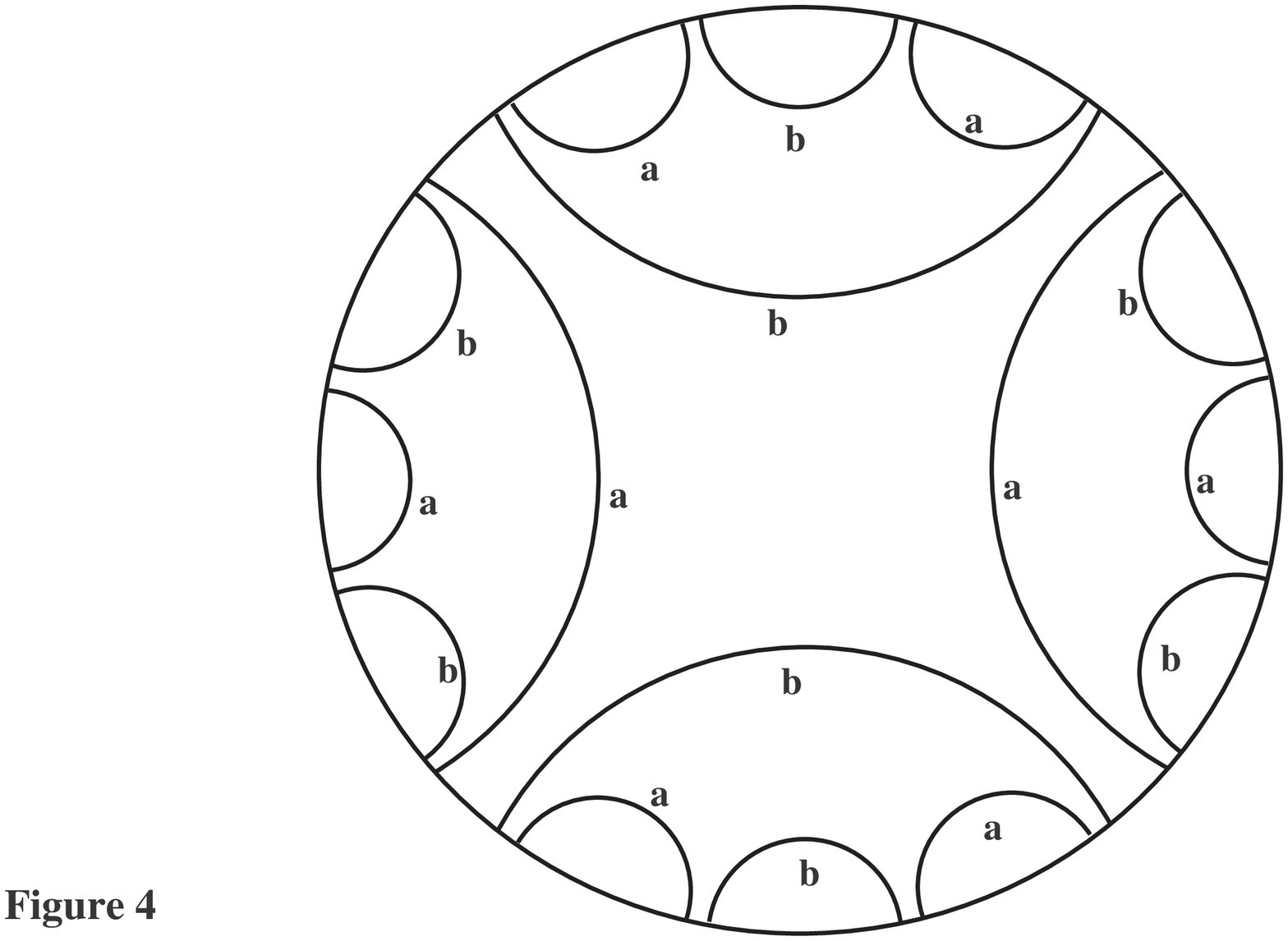,width=9cm,height=9cm,angle=-90}}

\end{figure}

Suppose that $\alpha$ is a geodesic segment or ray in $B^2$, which
does not lie in a translate of $a$ or $b$. Then $\a$ crosses a
sequence (finite or infinite, possibly of length $0$) of translates of
$a$ and $b$. When a geodesic segment starts or ends in a translate, or
a geodesic ray starts in one, that intersection is to be counted as a
crossing. To $\alpha$, we associate a sequence $S(\alpha)\eq
x_1x_2x_3\cdots$ of elements in the set $\set{a,\overline
{a},b,\overline {b}}$ in the following way. If the $i^{th}$ crossing
of $\alpha$ with the union of the translates of $a$ and $b$ is a
positive crossing with a translate of $a$, then $x_i\eq a$. If the
$i^{th}$ crossing is a negative crossing with a translate of $a$, then
$x_i\eq\overline {a}$. For crossings with translates of $b$, the
elements $b$ and $\overline {b}$ are assigned similarly. Note that
$S(\alpha)$ is an infinite sequence if and only if $\alpha$ is a
geodesic ray which ends at a limit point of $\Gamma_0$, and that for
each sequence $S\eq x_1x_2x_3\cdots$ of elements of the set
$\set{a,\overline {a},b,\overline {b}}$ with the property that for no
$i$ is $x_ix_{i+1}$ in the set
$\set{a\overline{a},\overline{a}a,b\overline{b},\overline{b}b}$, there
exist geodesic rays $\alpha$ with $S(\alpha)\eq S$.

Although we will not need it, the following fact seems worth
mentioning. If $S(\alpha)\eq x_1x_2x_3\ldots$ is a crossing sequence
of infinite length, let $\sigma(S(\alpha))$ denote the shifted
sequence $x_2x_3x_4\ldots$. Suppose $\a_1$ and $\alpha_2$ are geodesic
rays ending at the limit points $p_1$ and $p_2$ respectively. Then
there exist $m,n\geq 0$ so that
$\sigma^m(S(\alpha_1))\eq\sigma^n(S(\alpha_2))$ if and only if there
exists $\gamma\in\Gamma_0$ so that $\gamma(p_1)\eq p_2$.

Using these sequences, the controlled concentration points of $\Gamma_0$
can be characterized. The following result appears in \cite{A-H-M},
but we reproduce its proof here for the convenience of the reader.

\begin{proposition}{4.1}{}
Let $p$ be a limit point of $\Gamma_0$ which is the endpoint of a
geodesic ray $\alpha$ with $S(\alpha)\eq x_1x_2x_3\cdots$. Then $p$ is
a controlled concentration point for $\Gamma_0$ if and only if
$S(\alpha)$ has the following property. There exists $N$ such that for
all $n\geq N$, for all positive $k$, and for all $M$, there exists
$m\geq M$ such that $x_{n+i}\eq x_{m+i}$ for all $i$ with $0\leq i\leq
k$.
\marginwrite{4.1}
\end{proposition}

\noindent In words, past some point every finite subsequence reappears
infinitely many times. This is equivalent to the condition that past
some point every finite subsequence reappears at least once.

\begin{proof}{} Denote by $\lambda_n$ the translate of $a$ or $b$
whose crossing with $\alpha$ determines $x_n$, and by $U_n$ the
neighborhood of $p$ bounded by the endpoints of $\lambda_n$. Suppose
the condition in the proposition holds. By truncating $\alpha$, we may
assume that every subsequence reappears infinitely often. Let $m_k$ be
an integer greater than $k$ so that $x_{1+i}\eq x_{m_k+i}$ for $0\leq
i\leq k$. Given a neighborhood $V$ of $p$, choose $k$ so large that
$\lambda_k$ has endpoints in $V$. Let $\gamma$ be the element of
$\Gamma_0$ that translates $\lambda_1$ to $\lambda_{m_k}$. Note that
$\gamma$ translates $\lambda_{1+i}$ onto $\lambda_{m_k+i}$ for all
$0\leq i\leq k$, hence translates $U_{1+i}$ onto $U_{m_k+i}$ for
$0\leq i\leq k$. Therefore $\gamma(U_1)\eq U_{m_k}\subseteq V$ and
$p\in U_{m_k+k}\eq \gamma(U_{1+k})\subseteq \gamma(V)$, showing that
$U_1$ can be concentrated with control. Conversely, suppose $p$ is a
controlled concentration point and choose $N$ large enough so that
$U_N$, and hence every neighborhood of $p$ inside $U_N$, can be
concentrated with control. For any $n,k>N$ and any $M$, there exists
$\gamma$ so that $\gamma(U_n)\subseteq U_M$ and $\gamma^{-1}(p)\in
U_{n+k}$. This $\gamma$ must move
$\lambda_n,\lambda_{n+1},\ldots,\lambda_{n+k}$ onto a sequence of
translates of $a$ and $b$ crossed by $\alpha$, with endpoints in
$U_M$. Thus the condition of proposition~\ref{4.1} holds.
\end{proof}

Proposition~\ref{4.1} shows immediately that not all limit points of
$\Gamma_0$ are controlled concentration points. The next two theorems
provide more delicate examples.

\begin{theorem}{4.2}{}
There are uncountably many limit points of $\Gamma_0$ which are
concentration points but are not controlled concentration points.
\marginwrite{4.2}
\end{theorem}

\begin{proof}{} Denote by $a_n$ a sequence of $n$ $a$'s, and by
${\overline {a}}_n$ a sequence of $n$ \hskip0pt $\overline {a}$'s.
Choose one of the uncountably many increasing sequences of positive
integers $1\leq i_1<j_1<i_2<j_2<i_3<\cdots$, and let $p$ be a limit
point which is the endpoint of a geodesic ray whose crossing sequence
is
$$ba_{i_1}b\overline a_{j_1}ba_{i_2}b\overline a_{j_2}
ba_{i_3}b\overline a_{j_3}ba_{i_4}b\overline a_{j_4}\cdots\ .$$

\noindent By proposition~\ref{4.1}, $p$ is not a controlled
concentration point. We will verify that it is a concentration point.

Let $W$ be the connected neighborhood of $p$ in $\partial \overline
{B^2}$ whose endpoints are the translate of $b$ whose intersection
with $\alpha$ corresponds to the first $b$ in the crossing sequence of
$\alpha$. If $\lambda$ is any geodesic whose endpoints are the
endpoints of a connected neighborhood $U$ of $p$ with $U\subseteq W$,
then $\alpha$ crosses $\lambda$. We will show that any such $U$ can be
concentrated at~$p$.

We refer to Figure~5. The geodesic labeled with $a_{i_n}$ represents
the $i_n$ translates of $a$  whose positive crossings by $\alpha$
produce the block of $i_n$ $a$'s in the crossing sequence of $\a$, and
similarly for the geodesic labeled with $a_{j_n}$. The geodesic $\mu$
is the unique oriented geodesic which crosses the middle one of the
three translates of $b$ in Figure~5 which cross $\alpha$, and has
two-ended crossing sequence $\cdots a a a b \overline {a} \overline
{a} \overline {a} \cdots$~. Its endpoints are labeled $B$ and $C$, and
also shown are the two translates of $\mu$ which cross the other two
translates of $b$ that cross $\alpha$ in Figure~5. Since the crossing
sequence of $\alpha$ contains arbitrarily long blocks of the form
$a_{i_m}b\overline a_{j_m}$, it follows that there are translates of
$\alpha$ limiting onto $\mu$ so that the images of the initial point
of $\alpha$ limit to $B$ and the images of $p$ limit to $C$. This is
indicated by the direction arrow on~$\mu$. Similarly, directions are
labeled on the other two translates of~$\mu$.

\begin{figure}[htb]

\centerline{
\epsfig{figure=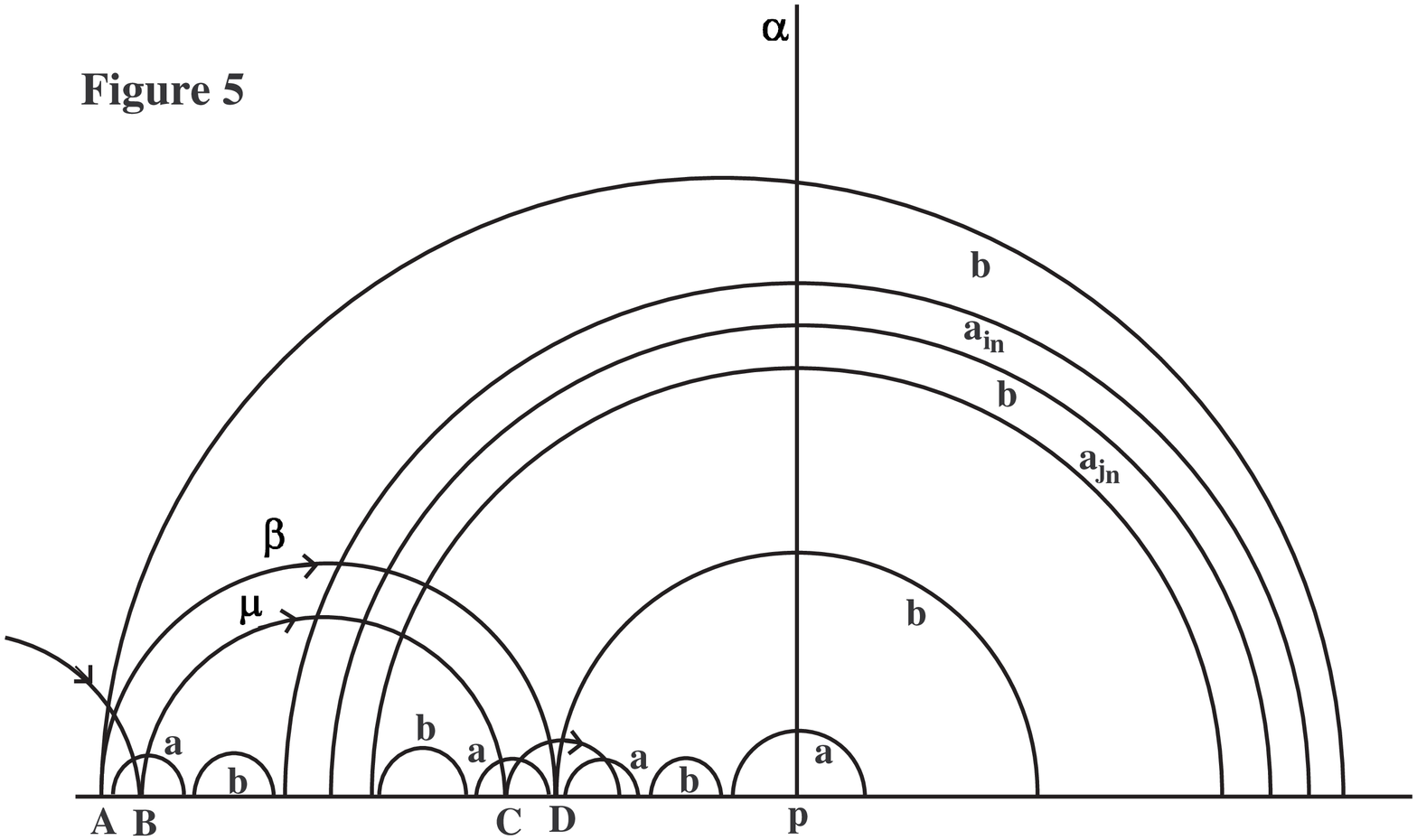,width=10cm,height=10.5cm,angle=-90}}
\end{figure}

We orient $\lambda$ so that it crosses $\alpha$ from left to right in
Figure~5. Let $E$ be the initial point of $\lambda$. There must be an
$n$ so that $E$ lies between the outermost and innermost of the three
translates of $b$ that cross $\a$ in Figure~5. Suppose first that $E$
is not equal to either $B$ or $C$. Then $\lambda$ must cross either
$\mu$ or one of the two translates of $\mu$ shown in Figure~5, making
some nonzero angle $\theta$ at the intersection point. Therefore it
crosses almost all the translates of $\alpha$ that limit to $\mu$.
Translating these back to $\alpha$, we find translates of $\lambda$
crossing $\alpha$ from left to right at angles approximately $\theta$,
arbitrarily close to $p$, showing that $U$ can be concentrated at~$p$.

There remains the case where $E$ equals either $B$ or $C$. For this,
consider the geodesic $\beta$ which runs from the left hand endpoint
$A$ of the largest translate of~$b$ in Figure~5 that crosses $\alpha$
to the left hand endpoint $D$ of the smallest one. Figure~6 shows the
images of $A$, $B$, $C$, and $D$ under the element $\gamma$ of
$\Gamma_0$ that moves the middle translate of $b$ in Figure~5 to a
corresponding one closer to $p$. To verify that the images of $A$ and
$D$ are as indicated, note that after crossing the translate of~$b$ in
Figure~5, $\beta$ makes $j_n$ negative crossings with $a$'s then
limits onto the unlabeled side of a translate of~$b$. Referring back
to the fundamental domain shown in Figure~4, one sees that the latter
translate of~$b$ must be labeled as shown in Figure~6, hence
$\gamma(D)$ is as shown. The determination of $\gamma(A)$ is similar.
So $\gamma(\lambda)$ runs from $\gamma(B)$ or $\gamma(C)$ to some
point lying between $\gamma(D)$ and $\gamma(A)$. This shows that $U$
can be concentrated at $p$, and completes the proof that $p$ is a
concentration point.
\end{proof}

\begin{figure}[htb]

\centerline{
\epsfig{figure=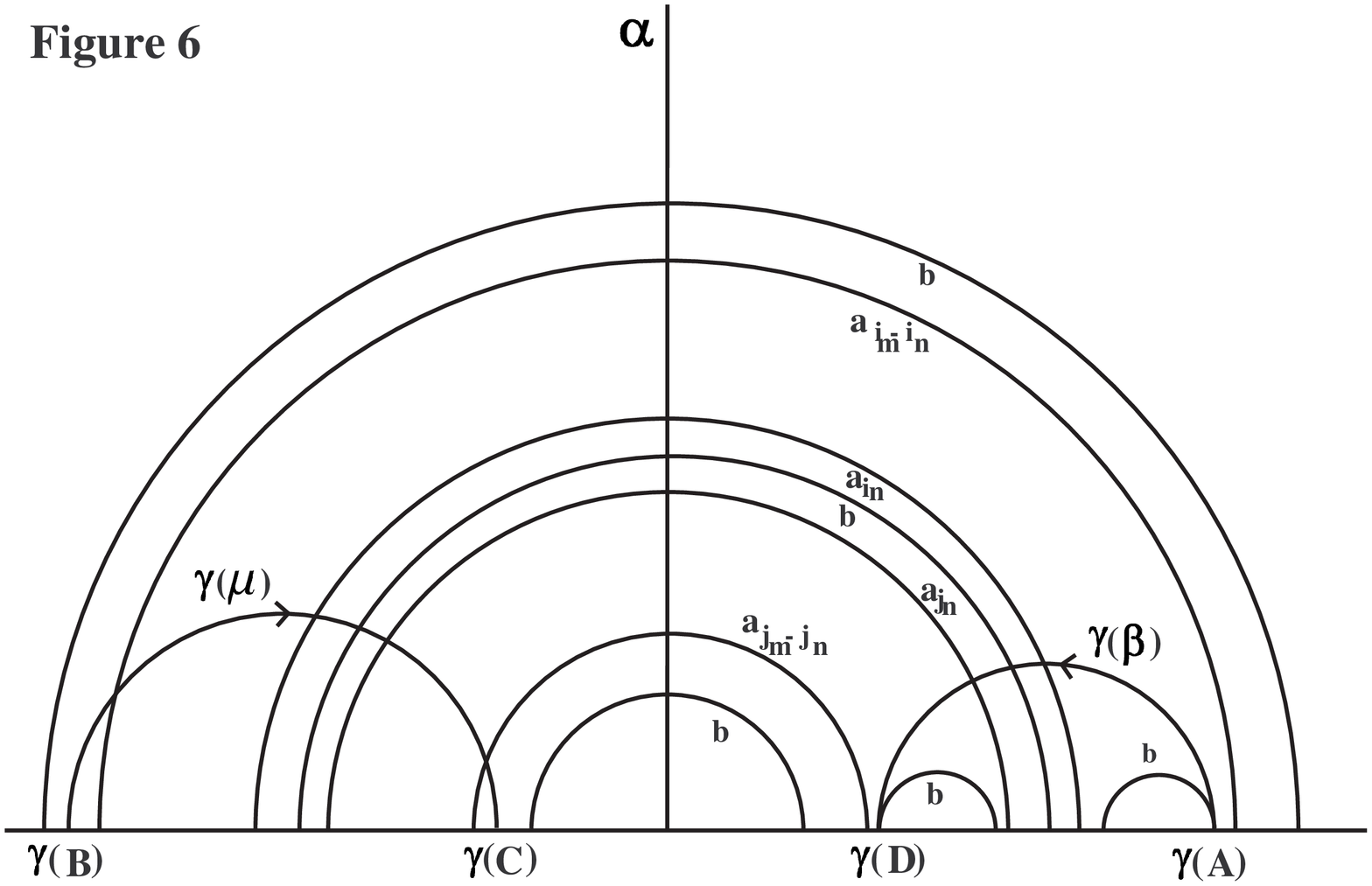,width=10cm,height=10.5cm,angle=-90}}
\end{figure}

\begin{theorem}{4.3}{}
There are uncountably many limit points of $\Gamma_0$ which are not
concentration points.
\marginwrite{4.3}
\end{theorem}

\begin{proof}{} We retain the notation of the proof of
theorem~\ref{4.2}. Choose one of the uncountably many increasing
sequences of positive integers $1\leq i_1<i_2<i_3<\cdots$, and let $p$
be a limit point which is the endpoint of a geodesic ray whose
associated sequence is $$ba_{i_1}ba_{i_2}ba_{i_3}ba_{i_4}\cdots\ .$$
\noindent We will verify that $p$ is not a concentration point.

We refer to Figure~7. Let $\lambda_n$ be the geodesic which runs from
the left hand endpoint of the larger translate of~$b$ in Figure~7(a)
to the right hand endpoint of the smaller one, and let $U_n$ be the
neighborhood of $p$ whose endpoints are the endpoints of $\lambda_n$.
We will show that $U_n$ cannot be concentrated at $p$. Since there are
arbitrarily small such neighborhoods, this implies that $p$ is not a
concentration point.

\begin{figure}[htb]

\centerline{
\epsfig{figure=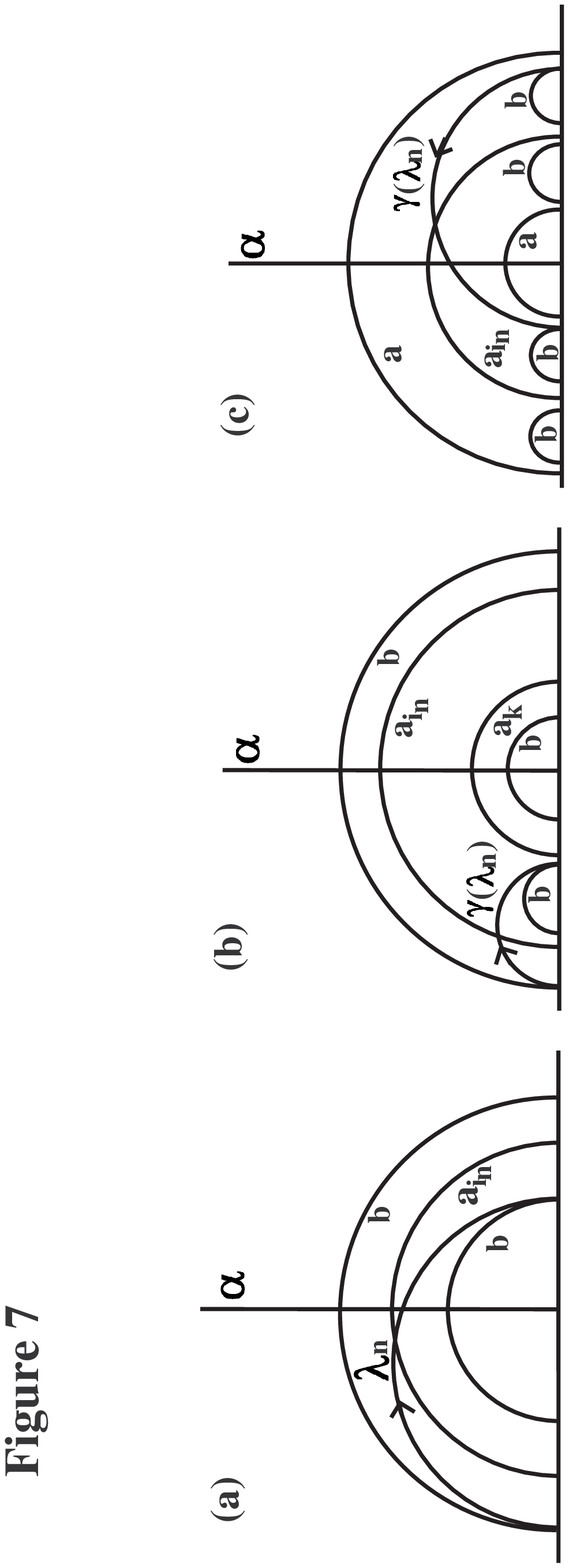,width=8cm,height=13.5cm,angle=-90}}

\end{figure}

Call the larger translate of~$b$ in Figure~7 $\mu_1$, and the smaller
$\mu_2$. Suppose there is an element $\gamma\in\Gamma_0$ so that
$\gamma(\lambda_n)$ crosses $\alpha$ from left to right, closer to $p$
(i.~e.~below the crossing of $\mu_2$ with $\alpha$). Suppose that
$\gamma(\mu_1)$ crosses $\alpha$. Then since $\lambda_n$ lies on the
unlabelled side of $\mu_1$, $\gamma(\mu_2)$ must lie underneath
$\gamma(\mu_1)$, and as shown in Figure~7(b), $\gamma(\lambda_n)$ must
lie entirely to the left of $p$. If $\gamma(\mu_2)$ crosses $\a$, then
similar considerations show that $\gamma(\lambda_n)$ lies to the right
of $p$. Suppose neither crosses $\alpha$. Observe that any translate
of $b$ which lies in the boundary of a translate of the fundamental
domain that intersects $\alpha$ near $p$ either
\begin{enumerate}
\item[\rm{(a)}] crosses $\alpha$, or
\item[\rm{(b)}] lies entirely to the left of $p$ and has its labelled side
underneath, or
\item[\rm{(c)}] lies entirely to the right of $p$ and has its labelled side
above.
\end{enumerate}

\noindent Since $\lambda_n$ crosses only $i_n$ translates of $\alpha$,
any translate of $\lambda_n$ that crosses $\alpha$ from left to right
near $p$ must either start on the unlabelled side of a translate of
$b$ or end on a labelled side of one, but neither of these is
possible. Note, however, that there are translates of $\lambda_n$
arbitrarily close to $p$ that cross $\alpha$ from right to left, as
shown in Figure~7(c). This is consistent with the fact that, by
theorem~\ref{3.1}, $p$ must be a geodesic separation point.
\end{proof}

\section [Another Example] {The infinitely generated case}
\label{another}

In this section we construct an infinitely generated Fuchsian group
$\Gamma_1$ having uncountably many conical limit points which are not
geodesic separation points. This shows that in theorem~\ref{3.1} the
hypothesis that $\Gamma$ is finitely generated is necessary. Moreover,
$\Gamma_1$ has uncountably many limit points that are weak
concentration points but not conical limit points.

\begin{proposition}{5.1}{}
There is an infinitely generated fuchsian group $\Gamma_1$, containing
no parabolic elements, having uncountably many weak concentration
points that are not conical limit points, and uncountably many conical
limit points that are not geodesic separation points.
\marginwrite{5.1}
\end{proposition}

\begin{proof}{} Let $\Gamma$ be the fundamental group of the
closed orientable surface $F$ of genus~2, acting on $B^2$ as
determined by some hyperbolic structure on $F$. It contains no
parabolic elements. Regard $F$ as the boundary of the genus~2
handlebody $V$, and choose elements $a$ and $b$ in $\pi_1(F)$ whose
images under the homomorphism $\pi_1(F)\to\pi_1(V)$ represent free
generators of $\pi_1(V)$.

Let $\widetilde V$ be the infinite cyclic covering of $V$
corresponding to the kernel of the homomorphism $\pi_1(V)\to\Z$ that
sends $a$ to $1$ and $b$ to $0$. This covering can be constructed by
cutting $V$ apart along a cocore disc $D_a$ for one of its 1-handles
(the one corresponding to the generator $a$) and gluing infinitely
many copies $\ldots,V_{-2}$, $V_{-1}$, $V_0$, $V_1$, $V_2,\ldots\,$,
of the split-open handlebody end to end along their copies of $D_a$.
For each $i$, $V_i\cap V_{i+1}$ is a lift $D_a^i$ of $D_a$. The cocore
disc $D_b$ for the other 1-handle lifts to a copy $D_b^i$ in each
$V_i$. Since every simple closed essential loop in $F$ is isotopic to
a geodesic, we may assume that $\partial D_a$ and $\partial D_b$ are
geodesics.

Let $\widetilde F$ be the boundary of $\widetilde V$, and let
$\Gamma_1$ be the subgroup of $\Gamma$ corresponding to
$\pi_1(\widetilde F)$. Denote $\widetilde F\cap V_i$ by $F_i$. Each
$F_i$ is a twice-punctured torus, with boundary $\partial
D_a^{i-1}\cup\partial D_a^i$, and with a 1-handle which contains the
loop~$\partial D_b^i$. Fix a basepoint $x$ of $F$, disjoint from
$\partial D_a$, and let $\widetilde x$ be the point of the preimage of
$x$ that lies in~$F_0$. Choose a basepoint $\widehat x$ for $B^2$ that
maps to $\widetilde x$. Notice that the union over all $j\in\Z$ of the
preimage geodesics of $\partial D_a^j\cup \partial D_b^j$ in $B^2$
forms the full preimage of $\partial D_a\cup \partial D_b$ in $B^2$.
Therefore these preimage geodesics are pairwise disjoint, and for
every $\epsilon>0$, there are only finitely many with Euclidean
diameter greater than~$\epsilon$.

Since $\widetilde F$ is a regular covering of $F$, $\Gamma_1$ is a
normal subgroup of $\Gamma$ and hence its limit set is all of $S^1$.
By corollary~3.3 of~\cite{M1}, every limit point of $\Gamma_1$ is a
weak concentration point. We will show that uncountably many of these
are not conical limit points.

For each $k$ let $c_k$ be the shortest loop based at $x$ that
represents $a^kb$ in $\pi_1(F,x)$. Choose uncountably many sequences
$i_1$, $i_2\ldots\,$ of positive integers, so that no two of the
sequences become equal after truncation of any initial segments. For
each sequence, let $\beta$ be the ray in $F$ corresponding to the
infinite product $c_{i_1}c_{i_2}\cdots\,$. Let $\alpha_0'$ be the lift
of $\beta$ to $\widetilde F$ starting at $\widetilde x$. For each
$j\geq 0$, $\alpha_0'$ crosses $\partial D_a^j$ exactly once.
Therefore the lift of $\alpha_0'$ to $B^2$, starting at $\hat x$,
limits to a single point $p$ in $S^1$. The geodesic ray $\alpha$ from
$\widehat x$ to $p$ also crosses the union of the preimage geodesics
of $\partial D_a^j$ exactly once. Let $\alpha_0$ be its image in
$\widetilde F$. Then for each $j\geq 0$, $\alpha_0$ crosses $\partial
D_a^j$ exactly once. In particular, for every compact subset $K$ of
$\widetilde F$, there is a subray of $\alpha_0$ which is disjoint from
$K$. This shows that $p$ is not a conical limit point for $\Gamma_1$.
The points $p$ obtained from different sequences are distinct. In
fact, no two of them can even be equivalent under the action of
$\Gamma_1$, for if so then some terminal segments of their defining
sequences would have to be equal.

To see that uncountably many limit points of $\Gamma_1$ are conical
limit points that are not geodesic separation points, modify the
previous construction by letting $c_k$ be the loop in $F$ representing
$a^{k}ba^{-k}b$. This time, $c_k$ lifts to a loop in $\widetilde F$
that starts at $\widetilde x$, moves into $F_k$, circles around the
1-handle in $F_k$ crossing $\partial D_b^k$, returns to $F_0$, and
goes around the 1-handle in $F_0$, crossing $\partial D_b^0$ once
before returning to $\widetilde x$. Choose the $i_k$ to be positive
and increasing, and proceed with the construction of $\beta$,
$\alpha_0'$, $\alpha$, and $\alpha_0$ as before. The fact that
$\alpha_0'$ crosses each $\partial D_b^{i_k}$ exactly once shows that
$\alpha_0'$ limits to a single point $p\in S^1$, and the geodesic ray
$\alpha_0$ crosses each $\partial D_b^{i_k}$ exactly once, since
$\alpha_0'$ does. Note that the crossing angles of $\alpha_0$ with the
$\partial D_b^{i_k}$ will be bounded away from~$0$. This time, $p$ is
a conical limit point, since $\alpha_0$ returns infinitely many times
to the compact subset $\partial D_b^0$. But $p$ is not a geodesic
separation point. For the crossing of $\alpha_0$ with $\partial
D_b^{i_k}$ produces a crossing of $\alpha$ with a geodesic $\lambda_k$
in the preimage of $\partial D_b^{i_k}$. Since the crossing angles are
bounded away from $0$, the endpoints of the $\lambda_k$ converge to
$p$. Since $\alpha_0$ crosses $\partial D_b^{i_k}$ only once, the
geodesics $\lambda_k$ show that $p$ is not a geodesic separation
point.
\end{proof}

\end{document}